\documentclass[12pt]{amsart}

\usepackage{amssymb}
\usepackage{amsmath}
\usepackage{amsfonts}
\usepackage{graphicx}
\usepackage{amsthm}
\usepackage{enumerate}
\usepackage[mathscr]{eucal}
\usepackage{mathrsfs}
\usepackage{verbatim}
\usepackage{booktabs}
\usepackage{enumitem}

\makeatletter
\@namedef{subjclassname@2020}{%
  \textup{2020} Mathematics Subject Classification}
\makeatother

\numberwithin{equation}{section}

\theoremstyle{plain}
\newtheorem{theorem}{Theorem}[section]

\newtheorem*{theorema}{Theorem A}
\newtheorem*{theoremb}{Theorem B}

\theoremstyle{remark}
\newtheorem{remark}[theorem]{Remark}

\DeclareMathOperator{\supp}{supp}

\renewcommand{\Re}{\operatorname{Re}}

\allowdisplaybreaks[4]

\begin{document}

\title[Endpoint estimates]{Some notes on endpoint estimates for pseudo-differential operators}


\author[J. Guo]{Jingwei Guo}
\address{Department of Mathematics\\
University of Science and Technology of China\\
Hefei 230026, China}
\email{jwguo@ustc.edu.cn}

\author[X. Zhu]{Xiangrong Zhu}
\address{Department of Mathematics\\
Zhejiang Normal University\\
Jinhua 321004, China}
\email{zxr@zjnu.cn}

\thanks{Xiangrong Zhu (the corresponding author) was supported by the NSFC Grant (No. 11871436). Jingwei Guo was supported by the NSF of Anhui Province, China (No. 2108085MA12).}

\date{\today}

\begin{abstract}
We study the pseudo-differential operator
 \begin{equation*}
 T_a f\left(x\right)=\int_{\mathbb{R}^n}e^{ix\cdot\xi}a\left(x,\xi\right)\widehat{f}\left(\xi\right)\,\textrm{d}\xi,
 \end{equation*}
 where the symbol $a$ is in the H\"{o}rmander class $S^{m}_{\rho,1}$ or more generally in the rough H\"{o}rmander class $L^{\infty}S^{m}_{\rho}$ with $m\in\mathbb{R}$ and $\rho\in [0,1]$. It is known that $T_a$ is bounded on $L^1(\mathbb{R}^n)$ for $m<n(\rho-1)$. In this paper we mainly investigate its boundedness properties when $m$ is equal to the critical index $n(\rho-1)$.

For any $0\leq \rho\leq 1$ we construct a symbol $a\in S^{n(\rho-1)}_{\rho,1}$ such that $T_a$ is unbounded on $L^1$ and furthermore it is not of weak type $(1,1)$ if $\rho=0$. On the other hand we prove that $T_a$ is bounded from $H^1$ to $L^1$ if $0\leq \rho<1$ and construct a symbol $a\in S^0_{1,1}$ such that $T_a$  is unbounded from $H^1$ to $L^1$.

Finally, as a complement, for any $1<p<\infty$ we give an example $a\in S^{-1/p}_{0,1}$ such that $T_a$ is unbounded on $L^p(\mathbb{R})$.
\end{abstract}



\maketitle

 \section{Introduction} \label{intro}

 Pseudo-differential operators are used extensively in the theory of partial differential equations and quantum field theory. The study of pseudo-differential operators began with the work of Kohn and Nirenberg \cite{KN65} and H\"{o}rmander \cite{H65}. A pseudo-differential operator is an operator given by
 \begin{equation*}
 T_a f(x)=\int_{\mathbb{R}^n}e^{ix\cdot\xi}a(x,\xi)\widehat{f}(\xi)\,\textrm{d}\xi, \qquad f\in \mathscr{S}(\mathbb{R}^n),
 \end{equation*}
where  $\widehat{f}$ is the Fourier transform of $f$ and  the symbol $a$ belongs to a certain symbol class. One of the most important symbol classes is the H\"{o}rmander class
 $S^m_{\rho,\delta}$ introduced in H\"{o}rmander \cite{H66}. A function $a\in C^{\infty}(\mathbb{R}^n\times\mathbb{R}^n)$ belongs to $S^{m}_{\rho,\delta}$ $(m\in \mathbb{R},0\leq\rho,\delta\leq1)$, if for all multi-indices $\alpha$ and $\beta$ it satisfies
 \begin{equation*}
 \sup_{x,\xi\in\mathbb{R}^{n}}(1+|\xi|)^{-m+\rho|\alpha|-\delta|\beta|}
 \left|\partial^{\alpha}_{\xi}\partial^{\beta}_{x}a(x,\xi)\right|<+\infty.
 \end{equation*}

 For pseudo-differential operators, one of the most important problems is whether they are bounded on the Lebesgue space $L^p(\mathbb{R}^n)$  and the Hardy space $H^1(\mathbb{R}^n)$. This problem has been extensively studied and there are numerous results. To mention a few, if $a\in S^{m}_{\rho,\delta}$ with $\delta<1$, $m\leq 0$ and $m\leq n(\rho-\delta)/2$, then $T_a$ is bounded on $L^2$. See H\"{o}rmander \cite{H71}, Hounie \cite{H86}, etc. For $a\in S^{m}_{\rho,1}$ Rodino \cite{R76} proved that $T_a$ is bounded on $L^2$ if  $m<n(\rho-1)/2$. He also constructed a symbol $a\in S^{m}_{\rho,1}$ such that $T_a$ is unbounded on $L^2$ if $m=n(\rho-1)/2$. (For the failure of $L^2$ boundedness for symbols in $S^{0}_{1,1}$, see also Ching \cite{Ching} and Stein \cite[P. 272]{S93}.) For endpoint estimates, in some unpublished lecture notes, Stein showed that if $a\in S^{n(\rho-1)/2}_{\rho,\delta}$ and either $0\leq\delta<\rho=1$ or  $0<\delta=\rho<1$, then $T_a$ is of weak type $(1, 1)$ and bounded from $H^1$ to $L^1$. This result was extended in \'{A}lvarez and Hounie \cite{AH90} to
\begin{theorema}[{\cite[Theorem 3.2]{AH90}}]
  If $a\in S^{m}_{\rho,\delta}$, $0<\rho\leq 1$, $0\leq \delta<1$ and $m=\frac n2(\rho-1+\min\{0,\rho-\delta\})$, then $T_a$ is of weak type $(1, 1)$ and bounded from $H^1$ to $L^1$.
 \end{theorema}
\noindent The $L^p$ boundedness then follows from an interpolation argument.

It is interesting to further investigate properties (say, $L^1$ boundedness) of $T_a$ for $a\in S^{m}_{\rho,1}$ (namely, $\delta=1$).  In fact one can consider more general symbol classes which contain $S^{m}_{\rho,1}$, for instance, the rough H{\"o}rmander class $L^{\infty}S^{m}_{\rho}$.

A systematic study of pseudo-differential operators with symbols with limited smoothness was initiated by Taylor in \cite{T91}.
 Rough symbols without any regularity assumption in the spatial variable were explored in Kenig and Staubach \cite{KW07} and Stefanov \cite{S10}.  In particular, Kenig and Staubach \cite{KW07} defined and studied  the following class of rough symbols which behave in the spatial variable $x$ like an $L^{\infty}$ function. A function $a(x,\xi)$, which is smooth in the frequency variable $\xi$ and bounded measurable in the spatial variable $x$, belongs to the rough H{\"o}rmander class $L^{\infty}S^{m}_{\rho}$ ($m\in \mathbb{R}$, $0\leq\rho\leq 1$), if for all multi-indices $\alpha$ it satisfies
 \begin{equation*}\sup_{\xi\in \mathbb{R}^n}(1+|\xi|)^{-m+\rho|\alpha|}\left\|\partial^{\alpha}_{\xi}a(\cdot,\xi)\right\|_{L^{\infty}( \mathbb{R}^n)}<+\infty.
 \end{equation*}
They proved the following theorem.
 \begin{theoremb}[{\cite[Proposition 2.3]{KW07}}]
 Let $a\in L^{\infty}S^{m}_{\rho}, 0\leq\rho\leq 1$.
 Assume that $m<n(\rho-1)/p$ and $1\leq p\leq 2$, then the operator $T_a$ is bounded on $L^p$.
 \end{theoremb}
\noindent Rodr\'{\i}guez-L\'{o}pez and Staubach \cite{RLS13} continued the investigation of \cite{KW07}, considered
certain rough symbol classes (e.g. $L^p S^{m}_{\rho}$) and studied both linear and multilinear pseudo-differential and Fourier integral operators.


 Our main results are concerning the behavior of $T_a$ with the symbol $a$ belonging to the critical class $L^{\infty}S^{n(\rho-1)}_{\rho}$.

 On one hand, by constructing a counterexample we show that the $L^1$ boundedness of $T_a$ may fail and Theorem B is sharp if $p=1$.
 \begin{theorem}\label{1.1}
 For any $0\leq \rho\leq 1$ there exists a symbol $a\in S^{n(\rho-1)}_{\rho,1}\subset L^{\infty}S^{n(\rho-1)}_{\rho}$ such that $T_a$ is unbounded on $L^1$ and furthermore if $\rho=0$ then $T_a$ is not of weak type $(1,1)$.
 \end{theorem}

On the other hand, we confirm the $H^1$--$L^1$ boundedness of $T_a$.
 \begin{theorem}\label{1.2}
 If $a\in L^\infty S^{n(\rho-1)}_{\rho}$, $0\leq \rho<1$, then $T_a$ is bounded from $H^1$ to $L^1$. If $\rho=1$, there exists a symbol $a\in S^{0}_{1,1}\subset L^\infty S^0_1$ such that $T_a$ is unbounded from $H^1$ to $L^1$.
 \end{theorem}

\begin{remark}
For the symbol $a\in S^{0}_{1,1}$ we construct, it is also true that $T_a$ is unbounded on $L^p$, $1<p<\infty$. See Subsection \ref{subsec3}.

Theorem \ref{1.2} extends the symbol class of Theorem A for the $H^1$--$L^1$ boundedness from $S^{m}_{\rho,\delta}$ to  $L^\infty S^{n(\rho-1)}_{\rho}$ if $0\leq\rho<1$ and $0\leq \delta \leq 1$. In particular, the case $\rho=0$ and the case $\delta=1$ were not included in Theorem A.

Our proofs for $0<\rho<1$ and $\rho=0$ rely on different methods. See Subsection \ref{subsec1} and \ref{subsec2}.  Since $T_a$ is unbounded on $L^2$ for some $a\in L^\infty S^{n(\rho-1)/2}_{\rho}$, we cannot follow \'{A}lvarez and Hounie's method in \cite{AH90}. Instead, we take advantage of properties of both the kernel and the Fourier transform of functions from Hardy spaces and adopt some argument from Seeger, Sogge and Stein \cite{SSS91} to obtain the desired boundedness.

Unlike the case $\rho=0$ (exhibited in Theorem \ref{1.1}), we were not able to determine whether $T_a$ is of weak type $(1,1)$ if $a\in L^\infty S^{n(\rho-1)}_{\rho}$, $0<\rho<1$. The difficulty might come from the fact that $T_a$ is unbounded on $L^2$ for some $a\in S^{n(\rho-1)/2}_{\rho,1}$.
\end{remark}

We also want to know whether $T_a$ is bounded on $L^p(\mathbb{R}^n)$ with $a\in L^\infty S^{n(\rho-1)/p}_{\rho}$. In this paper for simplicity we only consider the case $n=1$. We show that it is unbounded on $L^p(\mathbb{R})$ if $\rho=0$ and $1<p<\infty$.
 \begin{theorem}\label{1.4}
If $n=1$ and $1<p<\infty$, then there exists a symbol $a\in S^{-1/p}_{0,1}\subset L^\infty S^{-1/p}_{0}$ such that $T_a$ is unbounded on $L^p$.
 \end{theorem}

\begin{remark}
The main point of this theorem is the unboundedness for $1<p<2$. If $p=2$ Rodino \cite{R76} had constructed a symbol $a\in S^{n(\rho-1)/2}_{\rho,1}$ such that $T_a$ is unbounded on $L^2$. If $p>2$ the critical index of $m$ is $n(\rho-1)/2$ rather than $n(\rho-1)/p$.

It seems hard to generalize our counterexample  from the case $\rho=0$ to the case $\rho>0$.
\end{remark}

{\it Notations:}  Throughout this paper we will use the following notations. Let $B_r$ be the ball in $\mathbb{R}^n$ centered at the origin with radius $r$. We will use a standard partition of unity (which appears in the Littlewood-Paley theory for example).  There exist nonnegative functions $\Psi_{0}\in C^{\infty}_{c}(B_{3/4})$ with $\Psi_{0}\equiv 1$ on $B_{2/3}$ and  $\psi\in C^{\infty}_{c}(B_{3/2}\!\!\setminus \!B_{2/3})$ with $\psi\equiv1$ on $B_{4/3}\!\!\setminus \!B_{3/4}$ such that
 \begin{equation*}
 \Psi_{0}(\xi)+\sum^{\infty}_{j=0}\psi(2^{-j}\xi)=1, \quad\textrm{for all $\xi\in \mathbb{R}^n$}.
 \end{equation*}
Define two Schwartz functions $\Phi_0$ and $\varphi$ by $\Psi_{0}=\widehat{\Phi_0}$ and $\psi=\widehat{\varphi}$.

 For functions $f$ and $g$ with $g$ taking nonnegative real values,
$f\lesssim g$ means $|f|\leq Cg$ for some constant $C$. If $f$ is nonnegative, $f\gtrsim g$ means $g\lesssim f$. The notation $f\asymp g$ means that $f\lesssim g$ and $g\lesssim f$. Implicit constants shown up in this paper may depend on $n$, $\rho$ and $p$.


\section{Proof of Theorem \ref{1.1}}

If $\rho=1$ we will construct a symbol $a\in S^0_{1,1}$ in Subsection \ref{subsec3} such that $T_a$ is unbounded from $H^1$ to $L^1$, hence unbounded on $L^1$.

We assume that $0\leq \rho<1$ in what follows.  Define
 \begin{equation*}
 a(x,\xi)=\int_{\mathbb{R}^n}|u|^{-n}\left(1-\Psi_0(u)\right)\Phi_0\left(|u|^{-\rho}(\xi-u)\right)e^{-i x\cdot u}\,\textrm{d}u.
 \end{equation*}

We first verify that $a\in S^{n(\rho-1)}_{\rho,1}$. Since $\partial^\alpha \Phi_0\in\mathscr{S}$, for any $N>1$ and multi-indices $\alpha$ and $\beta$ we readily get
\begin{equation}\label{eq111}
\left|\partial^\alpha_\xi\partial^\beta_x a(x,\xi)\right|\lesssim_{\alpha,N}\int_{|u|>2/3}|u|^{|\beta|-\rho|\alpha|-n}\left(1+|u|^{-\rho}|\xi-u|\right)^{-N}\,\textrm{d}u.
\end{equation}
If $|\xi|\leq 1/3$, we choose an integer $N$ large enough such that $|\beta|-\rho|\alpha|-N(1-\rho)<0$ (which is possible since $\rho<1$). Then \eqref{eq111} yields
 \begin{align*}
 \left|\partial^\alpha_\xi \partial^\beta_x a(x,\xi)\right| &\lesssim_{\rho,\alpha,\beta}\int_{|u|>2/3}|u|^{|\beta|-\rho|\alpha|-N(1-\rho)-n}\,\textrm{d}u
 \lesssim_{n, \rho,\alpha,\beta} 1\\
 &\lesssim_{n, \rho,\alpha,\beta} (1+|\xi|)^{n(\rho-1)-\rho|\alpha|+|\beta|}.
 \end{align*}
If $|\xi|>1/3$, we fix an integer $N>\max\{n,|\alpha|,\frac{|\beta|-\rho|\alpha|}{1-\rho}\}$ and split the integral on the right side of \eqref{eq111} into three parts over $\{u: 2/3<|u|\leq|\xi|/2\}$, $\{u: |\xi|/2<|u|\leq 2|\xi| \}$ and $\{u: |u|> 2|\xi| \}$ respectively. The first part is
\begin{align*}
&\lesssim_{n, \rho,\alpha,\beta} |\xi|^{-N}\int_{2/3<|u|\leq|\xi|/2}|u|^{|\beta|+\rho(N-|\alpha|)-n}\,\textrm{d}u\\
&\lesssim_{n, \rho,\alpha,\beta} |\xi|^{|\beta|+\rho(N-|\alpha|)-N}\ln (1+|\xi|)\lesssim_{n, \rho,\alpha,\beta} (1+|\xi|)^{n(\rho-1)-\rho|\alpha|+|\beta|},
\end{align*}
where we need  the factor $\ln(1+|\xi|)$ only when $\rho=|\beta|=0$. The second part is
\begin{align*}
&\lesssim_{n, \rho,\alpha,\beta} |\xi|^{|\beta|-\rho|\alpha|-n}\int_{|\xi|/2<|u|\leq 2|\xi|}(1+|\xi|^{-\rho}|\xi-u|)^{-N}\,\textrm{d}u\\
&\lesssim_{n, \rho,\alpha,\beta} |\xi|^{|\beta|-\rho|\alpha|-n}\int_{|z|\leq 3|\xi|}(1+|\xi|^{-\rho}|z|)^{-N}\,\textrm{d}z\\
&\lesssim_{n, \rho,\alpha,\beta} (1+|\xi|)^{n(\rho-1)-\rho|\alpha|+|\beta|}.
\end{align*}
The third part is
\begin{align*}
&\lesssim_{n, \rho,\alpha,\beta} \int_{|u|> 2|\xi|}|u|^{|\beta|-\rho|\alpha|-N(1-\rho)-n}\,\textrm{d}u\lesssim_{n, \rho,\alpha,\beta} |\xi|^{|\beta|-\rho|\alpha|-N(1-\rho)}\\
&\lesssim_{n, \rho,\alpha,\beta} (1+|\xi|)^{n(\rho-1)-\rho|\alpha|+|\beta|}.
\end{align*}
To sum up, if $|\xi|>1/3$ we also have
\begin{equation*}
\left|\partial^\alpha_\xi \partial^\beta_x a(x,\xi)\right|\lesssim_{n, \rho,\alpha,\beta} (1+|\xi|)^{n(\rho-1)-\rho|\alpha|+|\beta|}.
 \end{equation*}
Therefore we have  $a\in S^{n(\rho-1)}_{\rho,1}$.

We next show that $T_a$ is unbounded on $L^1$ if $0\leq \rho< 1$ and not of weak type $(1,1)$ if $\rho=0$. By using the definition of $a(x,\xi)$, changing variables and $\Psi_{0}=\widehat{\Phi_0}$, we can write
 \begin{equation*}
 T_a f(x)=\int_{\mathbb{R}^n}k(x,y)f(y)\,\textrm{d}y,
 \end{equation*}
 where
 \begin{equation*}
 k(x,y)=\int_{\mathbb{R}^n}e^{-iy\cdot u}|u|^{n(\rho-1)}(1-\Psi_0(u))\Psi_0(|u|^{\rho}(y-x))\,\textrm{d}u.
 \end{equation*}

If $0<\rho<1$,  let us consider $|y|<\epsilon$ and $3\epsilon^{\rho}<|x|<1/2$ for small $\epsilon$. Then
 \begin{equation*}
\Re k(x,y)\gtrsim\int_{1<|u|<(3|x-y|/2)^{-1/\rho}\leq \epsilon^{-1}}|u|^{n(\rho-1)}\,\textrm{d}u\gtrsim |x|^{-n}.
 \end{equation*}
 Set $f_\epsilon(y)=\epsilon^{-n}\chi_{B_1}(\epsilon^{-1}y)$. Then $\|f_\epsilon\|_1=|B_1|$. If $3\epsilon^{\rho}<|x|<1/2$ then  $|T_a f_\epsilon(x)|\gtrsim |x|^{-n}$ which implies that $\|T_a f_\epsilon\|_1\gtrsim \ln \epsilon^{-1}$. Thus $T_a$ is unbounded on $L^1$.

If $\rho=0$, let us consider $0<|y|<1/6$ and $|x|<1/2$. Then
 \begin{align*}
 \Re k(x,y)&=\int_{\mathbb{R}^n}\cos(y\cdot u)|u|^{-n}(1-\Psi_0(u))\,\textrm{d}u\\
 &\geq \int_{|u|>1}\cos(y\cdot u)|u|^{-n}\,\textrm{d}u-\int_{2/3<|u|\leq 1}|u|^{-n}\,\textrm{d}u
 \end{align*}
The second integral is an absolute constant. As to the first integral, by using polar coordinates, changing variables and splitting it into two, we have
 \begin{align*}
  &\int_{|u|>1}\cos(y\cdot u)|u|^{-n}\,\textrm{d}u=\int_{S^{n-1}}\int^\infty_1\cos(r\left|y\cdot \theta\right|)\frac{\textrm{d}r}{r}\textrm{d}\theta\\
 =&\int_{S^{n-1}}\int^{1}_{|y\cdot \theta|}\frac{\cos r}{r}\,\textrm{d}r\textrm{d}\theta+\int_{S^{n-1}}\int^\infty_1\frac{\cos r}{r}\,\textrm{d}r\textrm{d}\theta\\
 \geq &\cos 1\left| S^{n-1}\right|\ln |y|^{-1}+\left| S^{n-1}\right|\int^\infty_1\frac{\cos r}{r}\,\textrm{d}r,
 \end{align*}
where the last integral is an absolute constant as well. Therefore for small $|y|$
\begin{equation*}
 \Re k(x,y)\gtrsim \ln |y|^{-1}.
\end{equation*}
Consider $f_\epsilon(y)$ again. Then
 $|T_a f_\epsilon(x)|\gtrsim \ln \epsilon^{-1}$ for all $|x|<1/2$ and small $\epsilon$. Thus $\|T_a f_\epsilon\|_{L^{1,\infty}}\gtrsim\ln \epsilon^{-1}$ which  implies that $T_a$ is not of weak type $(1,1)$.  \qed

 \section{Proof of Theorem \ref{1.2}}

\subsection{Proof of Case $0<\rho<1$} \label{subsec1}
 It suffices to show $\|T_a b\|_1\lesssim 1$ for any $L^2$-atom $b$ for $H^1(\mathbb{R}^n)$ satisfying $\supp b\subset B_r$, $\int_{B_r}b(y)\,\textrm{d}y=0$ and $\|b\|_2\leq r^{-n/2}$.

As $a\in L^\infty S^{n(\rho-1)}_{\rho}$ we have a fundamental estimate of the kernel
\begin{equation}
|k(x,y)|=\left|\int_{\mathbb{R}^n}e^{i(x-y)\cdot \xi}a(x,\xi)\,\textrm{d}\xi\right|\lesssim_N |x-y|^{-N}\label{fd1}
\end{equation}
for any integer $N>n$. Indeed,  by integration by parts we have
 \begin{align*}
 |k(x,y)|&\lesssim |x-y|^{-N}\int_{\mathbb{R}^n}\left|\nabla^{N}_{\xi}a(x,\xi)\right|  \,\textrm{d}\xi \\
 &\lesssim |x-y|^{-N}\int_{\mathbb{R}^n}(1+|\xi|)^{n(\rho-1)-N\rho}\,\textrm{d}\xi\lesssim |x-y|^{-N}.
 \end{align*}
As a consequence, if $|x|>2r$ then
\begin{equation*}
 \left|T_a b(x)\right|=\left|\int_{|y|\leq r}k(x,y)b(y)\,\textrm{d}y\right|\lesssim_N\!\!\int_{|y|\leq r}|x-y|^{-N}|b(y)| \,\textrm{d}y\lesssim_N \!\!|x|^{-N}.
 \end{equation*}

If $r\geq1$, we take $N=n+1$ and get
 \begin{align*}
 \left\|T_a b\right\|_1&=\int_{|x|>2r}\left|T_a b(x)\right|\,\textrm{d}x+\int_{|x|\leq 2r}\left|T_a b(x)\right|\,\textrm{d}x\\
 &\lesssim \int_{|x|>2r}|x|^{-N}\,\textrm{d}x+r^{n/2}\left\|T_ab\right\|_2\\
 &\lesssim 1+r^{n/2}\|b\|_2\lesssim 1,
 \end{align*}
where we have used the $L^2$ boundedness of $T_a$ (ensured by \cite[Proposition 2.3]{KW07}).

We assume $r<1$ below. By using the partition of unity introduced in Section \ref{intro}, we divide $T_{a}b$ as
 \begin{align}
 T_{a}b(x)&=\int e^{ix\cdot\xi}a(x,\xi)\Psi_{0}(\xi)\widehat{b}(\xi)\textrm{d}\xi
 +\sum^{\infty}_{j=0}\int e^{ix\cdot\xi}a(x,\xi)\psi(2^{-j}\xi)\widehat{b}(\xi)\textrm{d}\xi\nonumber\\
 &=:T_{a\Psi_{0}}b(x)+\sum^{\infty}_{j=0}T_{j}b(x). \label{333}
 \end{align}
Notice that $T_{a\Psi_{0}}$ is just the pseudo-differential operator with symbol $a(x,\xi)\Psi_{0}(\xi)$. Since $\Psi_{0}\in\mathscr{S}$ we can use the $L^1$ boundedness of $T_{a\Psi_{0}}$ (ensured by \cite[Proposition 2.3]{KW07}) to obtain
\begin{equation*}
\|T_{a\Psi_{0}}b\|_1 \lesssim \|b\|_1\lesssim r^{n/2}\|b\|_2\leq 1.
\end{equation*}
Applying the $L^2$ boundedness of $T_a$ gives
\begin{equation*}
\int_{|x|\leq 3r}\left|\sum^{\infty}_{j=0}T_{j}b(x)\right|\,\textrm{d}x\lesssim 1.
\end{equation*}
Hence it suffices to prove
 \begin{equation}
 \sum^\infty_{j=0}\int_{|x|>3r}\left|T_{j}b(x)\right|\,\textrm{d}x\lesssim 1. \label{222}
 \end{equation}

If $2^j\geq r^{-1/\rho}$, we write
\begin{equation*}
T_{j}b(x)=\int_{\mathbb{R}^n}k_j(x,y)b(y)\,\textrm{d}y,
\end{equation*}
where
\begin{equation*}
k_j(x,y)=\int_{\mathbb{R}^n}e^{i(x-y)\cdot\xi}a(x,\xi)\psi(2^{-j}\xi)\,\textrm{d}\xi.
\end{equation*}
Set $L_{\xi}=1-2^{2j\rho}\triangle_\xi$. For an integer $N>n/2$ we have
 \begin{align*}
 \left|k_j(x,y)\right|&=\left(1+2^{2j\rho}|x-y|^2\right)^{-N}\left|\int_{\mathbb{R}^n}L_{\xi}^N\left(e^{i(x-y)\cdot\xi}\right) a(x,\xi)\psi(2^{-j}\xi)\,\textrm{d}\xi\right|\\
 &\leq \left(1+2^{2j\rho}|x-y|^2\right)^{-N}\!\!\int_{2^{j-1}<|\xi|<2^{j+1}}\left|L_{\xi}^N\left(a(x,\xi)\psi(2^{-j}\xi)\right)\right|\textrm{d}\xi.
 \end{align*}
Note that
 \begin{equation*}
 \left|L_{\xi}^N(a(x,\xi)\psi(2^{-j}\xi))\right|\lesssim 2^{jn(\rho-1)}
 \end{equation*}
which can be verified easily. Thus
 \begin{equation*}
 \left|k_j(x,y)\right|\lesssim (1+2^{2j\rho}|x-y|^2)^{-N}2^{jn\rho}.
 \end{equation*}
If $|y|<r$ and $|x|>3r$, then
\begin{equation*}
\left|T_jb(x)\right|\lesssim 2^{jn\rho}(1+2^{2j\rho}|x|^2)^{-N}.
\end{equation*}
Therefore
 \begin{equation*}
 \sum_{2^j\geq r^{-1/\rho}}\int_{|x|>3r}\left|T_{j}b(x)\right|\,\textrm{d}x
 \lesssim \sum_{2^j\geq r^{-1/\rho}}2^{j(n-2N)\rho}\int_{|x|>3r}|x|^{-2N}\,\textrm{d}x\lesssim 1.
 \end{equation*}

If $2^j<r^{-1/\rho}$, we first notice two basic estimates of the atom $b$: for any nonnegative integer $N$,
 \begin{equation}\label{H11}
 \left|\nabla^N\widehat{b}(\xi)\right|\leq \int_{|y|\leq r}|y|^N|b(y)|\,\textrm{d}y\lesssim_N r^N;
 \end{equation}
also, \cite[Section III.5.4]{S93} gives
 \begin{equation}\label{H12}
 \int_{\mathbb{R}^n}\big|\widehat{b}(\xi)\big||\xi|^{-n}\,\textrm{d}\xi\lesssim \|b\|_{H^1}\lesssim 1.
 \end{equation}
With $L_{\xi}$ defined as above and the definition of $T_{j}b$ in \eqref{333},  for any fixed integer $N>n/2$, we have
 \begin{align*}
\left|T_{j}b(x)\right|
&=(1+2^{2j\rho}|x|^2)^{-N}\left|\int L_{\xi}^N\left(e^{ix\cdot\xi}\right)a(x,\xi)\psi(2^{-j}\xi)\widehat{b}(\xi)\,\textrm{d}\xi\right|\\
 &\leq (1+2^{2j\rho}|x|^2)^{-N}\int \left|L_{\xi}^N\left(a(x,\xi)\psi(2^{-j}\xi)\widehat{b}(\xi)\right)\right|\,\textrm{d}\xi.
 \end{align*}
By the Leibniz rule and \eqref{H11},
\begin{align*}
         &\left|L_{\xi}^N\left(a(x,\xi)\psi(2^{-j}\xi)\widehat{b}(\xi)\right)\right|\\
\lesssim &\sum_{N_1+N_2+N_3\leq 2N}\left|2^{jN_1\rho}\nabla^{N_1}_\xi a(x,\xi)\right|\left|2^{jN_2\rho}\nabla_{\xi}^{N_2}\left(\psi(2^{-j}\xi)\right)\right|\left|2^{jN_3\rho}\nabla_{\xi}^{N_3}\widehat{b}(\xi)\right|\\
\lesssim & 2^{jn(\rho-1)} \sum_{0\leq N_3\leq 2N} 2^{jN_3\rho}\left|\nabla^{N_3}_\xi \widehat{b}(\xi)\right|\\
\lesssim & 2^{jn(\rho-1)} \left(\big|\widehat{b}(\xi)\big|+\sum_{1\leq N_3\leq 2N}\left(2^{j\rho}r\right)^{N_3}\right)\\
\lesssim & 2^{jn(\rho-1)} \left(\big|\widehat{b}(\xi)\big|+2^{j\rho}r\right).
\end{align*}
Thus, by \eqref{H12} we have
\begin{align*}
 &\quad \sum_{2^j<r^{-1/\rho}}\int_{\mathbb{R}^n}\left|T_{j}b(x)\right|\,\textrm{d}x\\
 &\lesssim \sum_{2^j< r^{-1/\rho}}2^{-jn}\int_{2^{j-1}<|\xi|<2^{j+1}}\left(\big|\widehat{b}(\xi)\big|+2^{j\rho}r\right)\,\textrm{d}\xi\\
 &\lesssim \sum_{2^j< r^{-1/\rho}}\left(2^{j\rho}r+\int_{2^{j-1}<|\xi|<2^{j+1}}\big|\widehat{b}(\xi)\big||\xi|^{-n}\,\textrm{d}\xi\right)\\
 &\lesssim  1+\int \big|\widehat{b}(\xi)\big||\xi|^{-n}\,\textrm{d}\xi \lesssim  1.
 \end{align*}
 Thus \eqref{222} is proved, as desired. \qed

 \subsection{Proof of Case $\rho=0$}\label{subsec2}
 In this part we show that $T_a$ is bounded from $H^1$ to $L^1$ if $a\in L^\infty S^{-n}_0$. Once again it suffices to show $\|T_a b\|_1\lesssim 1$ for any $L^2$-atom $b$ for $H^1(\mathbb{R}^n)$ satisfying $\supp b\subset B_r$, $\int_{B_r}b(y)\,\textrm{d}y=0$ and $\|b\|_2\leq r^{-n/2}$.

The method we use for this case is different from that for the previous case since the kernel is not rapidly decreasing any more.

If $r\geq 1$, we first use the $L^2$ boundedness of $T_a$ (ensured by \cite[Proposition 2.3]{KW07}) to get
\begin{equation*}
 \left\|T_a b\right\|_{L^1(B_{2r})}\lesssim r^{n/2}\left\|T_ab\right\|_2\lesssim r^{n/2}\|b\|_2\lesssim 1.
\end{equation*}
It remains to estimate the part over $\{x: |x|\geq 2r\}$. By integration by parts we get
\begin{equation*}
\left|T_a b(x)\right|=\left|x_l\right|^{-n}\left|\int_{\mathbb{R}^n} e^{ix\cdot\xi}\partial^{n}_{\xi_l}\left(a(x,\xi)\widehat{b}(\xi)\right)\,\textrm{d}\xi\right|
\end{equation*}
for some $1\leq l\leq n$ with $|x_l|\geq |x|/n$. By using this lower bound of $|x_l|$, a summation over $l$ and the Leibniz rule, for any nonzero $x$ we obtain
 \begin{align*}
 \left|T_ab(x)\right|&\lesssim |x|^{-n}\sum^{n}_{l=1}\sum^{n}_{j=0}\left|\int_{\mathbb{R}^n} e^{ix\cdot\xi}\partial^{n-j}_{\xi_l}a(x,\xi) \widehat{y^j_l b}(\xi)\,\textrm{d}\xi\right|\\
 &=:|x|^{-n}\sum^{n}_{l=1}\sum^{n}_{j=0}\left|T_{l,j}\left(y^j_l b\right)(x)\right|,
 \end{align*}
 where $T_{l,j}$ is the pseudo-differential operator with the symbol $\partial^{n-j}_{\xi_l}a$. It is obvious that $\partial^{n-j}_{\xi_l}a\in L^\infty S^{-n}_0$. Thus $T_{l,j}$ is bounded on $L^2$ and
 \begin{align*}
 \int_{|x|\geq 2r} \left|T_a b(x)\right|\,\textrm{d}x
 &\lesssim \sum^{n}_{l=1}\sum^{n}_{j=0}\int_{|x|\geq 2r} |x|^{-n}\left|T_{l,j}\left(y^j_l b\right)(x)\right|\,\textrm{d}x. \\
 &\lesssim \sum^{n}_{l=1}\sum^{n}_{j=0}r^{-n/2}\left\| y^j_l b(y)\right\|_{L^2(\textrm{d}y)} \\
 &\lesssim \sum^{n}_{j=0}r^{-n+j}\lesssim  1,
 \end{align*}
where we have used the assumption $r\geq1$. Therefore, we have proved $\|T_a b\|_1\lesssim 1$ if $r\geq1$.

If $r<1$,  we first use \eqref{H12} to get
\begin{equation}\label{eq4.3}
 \left|T_a b(x)\right|=\left|\int_{\mathbb{R}^n} e^{ix\cdot\xi}a(x,\xi)\widehat{b}(\xi)\,\textrm{d}\xi\right|\leq \int_{\mathbb{R}^n} \big|\widehat{b}(\xi)\big||\xi|^{-n}\,\textrm{d}\xi\lesssim 1.
 \end{equation}
It follows that
\begin{equation*}
\int_{|x|\leq 2} \left|T_a b(x)\right|\,\textrm{d}x\lesssim 1.
\end{equation*}
It remains to estimate the part over $\{x: |x|>2\}$. Arguing as above readily yields
\begin{align*}
 \left|T_ab(x)\right|&\lesssim |x|^{-2n}\sum^{n}_{l=1}\sum^{2n}_{j=0}\left|\int_{\mathbb{R}^n} e^{ix\cdot\xi}\partial^{2n-j}_{\xi_l}a(x,\xi) \widehat{y^j_l b}(\xi)\,\textrm{d}\xi\right|\\
 &=:|x|^{-2n}\sum^{n}_{l=1}\sum^{2n}_{j=0}\left|\widetilde{T}_{l,j}\left(y^j_l b\right)(x)\right|,
 \end{align*}
 where $\widetilde{T}_{l,j}$ is the pseudo-differential operator with the symbol $\partial^{2n-j}_{\xi_l}a$. Since  $\partial^{2n-j}_{\xi_l}a$ belongs to $L^\infty S^{-n}_0$,  we have that $\widetilde{T}_{l,j}$ is bounded from $L^{n/(n-1)}$ to $L^{n/(n-1)}$ if $n\geq 2$ and from $L^2$ to $L^2$ if $n=1$ (by \cite[Proposition 2.3]{KW07}). If $n\geq 2$ by using \eqref{H12}, \eqref{eq4.3} (with $a$ replaced by $\partial^{2n}_{\xi_l}a$) and H\"{o}lder's inequality we have
 \begin{align*}
 \int_{|x|>2} \left|T_ab(x)\right|\,\textrm{d}x&\lesssim \sum^{n}_{l=1}\sum^{2n}_{j=0}\int_{|x|>2}|x|^{-2n}\left|\widetilde{T}_{l,j}\left(y^j_l b\right)(x)\right|\,\textrm{d}x \\
 &\lesssim 1+\sum^{n}_{l=1}\sum^{2n}_{j=1}\left\| y^j_l b(y)\right\|_{L^{\frac{n}{n-1}}(\textrm{d}y)}\\
 &\lesssim 1+\sum^{2n}_{j=1}r^{j+\frac{n}{2}-1}\|b\|_2\lesssim  1.
 \end{align*}
If $n=1$ we just need to slightly modify the above computation to obtain the same bound by using the $L^2$ boundedness of $\widetilde{T}_{l,j}$.  Therefore, we have proved $\|T_a b\|_1\lesssim 1$ if $r<1$. This finishes the proof. \qed

 \subsection{Counterexample for $\rho=1$} \label{subsec3}
 We construct a counterexample $a\in S^0_{1,1}$ (similar to the example considered in Rodino \cite{R76})   such that $T_a$ is unbounded from $H^1$ to $L^1$. The existence of such an example is not surprising to us since $a\in S^0_{1,1}$ does not necessarily imply the $L^2$ boundedness of $T_a$.

Set $\theta=(1,0,\ldots,0)$ and
 \begin{equation*}
 a(x,\xi)=\sum^\infty_{j=2}e^{-i2^jx_1}\Psi_0\left(2^{-j+1}(\xi-2^j\theta)\right).
 \end{equation*}
It is routine to verify that $a\in S^0_{1,1}$. For any integer $N\geq2$, take $f_N$ such that
 \begin{equation*}
 \widehat{f_N}(\xi)=\sum^{N+1}_{s=2}\Psi_0\left(\xi-2^s\theta\right).
 \end{equation*}

We first show that
\begin{equation}
 \left\|T_a f_N\right\|_1\asymp N.\label{eq5.1}
\end{equation}
Indeed,  we notice that for any integer $j,s\geq 2$
 \begin{equation*}
 \Psi_0\left(2^{-j+1}(\xi-2^j\theta)\right)\Psi_0(\xi-2^s\theta)=\delta^s_j\Psi_0(\xi-2^s\theta),
 \end{equation*}
where $\delta^s_j$ is the Kronecker notation. Then
 \begin{equation*}
 T_af_N(x)=\sum^{N+1}_{s=2}\int_{\mathbb{R}^n} e^{ix\cdot(\xi-2^s\theta)}\Psi_0(\xi-2^s\theta)\,\textrm{d}\xi=N\widehat{\Psi_0}(-x),
 \end{equation*}
which immediately leads to \eqref{eq5.1}.

We next estimate the size of  $\|f_N\|_{H^1}$. Denote $\psi_j(\xi)=\psi(2^{-j}\xi)$ and $\varphi_j(x)=2^{jn}\varphi(2^{j}x)$ for any $j\in\mathbb{Z}$, where $\psi$ and $\varphi$ are as mentioned in Section \ref{intro}. By using the equivalence of norms of $H^p$ and $\dot{F}^{0,2}_p$ for $0<p\leq 1$ (see \cite[Section 6.5]{G09}) we have
 \begin{equation*}
 \|f_N\|_{H^1} \asymp \|f_N\|_{\dot{F}^{0,2}_{1}}=\left\|\left(\sum_{j\in\mathbb{Z}}\left|\varphi_j\ast f_N\right|^2\right)^{1/2}\right\|_{1}.
 \end{equation*}
Since it is easy to check that
 \begin{equation*}
 \psi_j(\xi)\Psi_0(\xi-2^s\theta)=\delta^s_j\Psi_0(\xi-2^s\theta)
 \end{equation*}
for $s\geq 2$ and $j\in\mathbb{Z}$,  we have
 \begin{equation*}
 \widehat{\varphi_j\ast f_N}(\xi)=\sum^{N+1}_{s=2}\psi_j(\xi)\Psi_0(\xi-2^s\theta)=\sum^{N+1}_{s=2}\delta^s_j\Psi_0(\xi-2^s\theta),
 \end{equation*}
 which yields that
 \begin{equation*}
 \varphi_j\ast f_N(x)=\left\{
                           \begin{array}{ll}
                           e^{i 2^j x_1}\Phi_0(x), &   \textrm{if $2\leq j \leq N+1$,}\\
                           0,                      &   \textrm{otherwise.}
                           \end{array}
                    \right.
 \end{equation*}
It follows that
 \begin{equation}
 \left\|f_N\right\|_{H^1}\asymp \left\|\left(\sum^{N+1}_{j=2}\left|e^{i2^jx_1}\Phi_0(x)\right|^2\right)^{\!\!1/2}\right\|_{L^1(\textrm{d}x)}\!\!\!\!\!
 =\left\|N^{1/2}\Phi_0\right\|_{1}\asymp N^{1/2}.\label{eq5.3}
 \end{equation}

 By \eqref{eq5.1} and \eqref{eq5.3}, for any integer $N\geq 2$ we have
 \begin{equation*}
 \left\|T_af_N\right\|_1\asymp N^{1/2}\left\|f_N\right\|_{H^1},
 \end{equation*}
 which implies that $T_a$ is unbounded from $H^1$ to $L^1$.   \qed

\begin{remark}
In fact, arguing as above easily yields that $T_a$ is unbounded on $L^p$ for any $1<p<\infty$ if we notice the equivalence of norms of $L^p$ and $\dot{F}^{0,2}_p$ for $1<p<\infty$ (see \cite[Section 6.5]{G09}).
 \end{remark}

 \section{Proof of Theorem \ref{1.4}}

For $x, \xi\in \mathbb{R}$, set
 \begin{equation*}
 a(x,\xi)=\sum^\infty_{j=2}2^{-\frac{j}{p}}\sum_{\frac{7}{8}2^{j}<|k|<\frac{9}{8}2^{j}} e^{-i4kx}\Psi_0\left(4^{-1}(\xi-4k)\right),
 \end{equation*}
where $k$ is implicitly assumed to be in $\mathbb{Z}$. We will follow this convention below.  It is easy to check that $a\in S^{-1/p}_{0,1}$. For any integer $N>8$, take $f_N$ such that
 \begin{equation*}
 \widehat{f_N}(\xi)=\sum^{N+1}_{s=2}2^{-s(1-\frac 1p)}\sum_{\frac{7}{8}2^{s}<|k|<\frac{9}{8}2^{s}}\Psi_0(\xi-4k).
 \end{equation*}

We first show that if $p<\infty$ then
 \begin{equation}\label{eq6.1}
 \left\|T_af_N\right\|_{p}\gtrsim N.
 \end{equation}
Indeed, since it is easy to check that
\begin{equation*}
 \Psi_0\left(4^{-1}(\xi-4k_1)\right)\Psi_0(\xi-4k_2)=\delta^{k_1}_{k_2}\Psi_0(\xi-4k_2)
\end{equation*}
for $k_1,k_2\in \mathbb{Z}$, where $\delta^{k_1}_{k_2}$ is the Kronecker notation,  we then get
 \begin{equation*}
    \left|T_af_N(x)\right|=\left|\Phi_0(x)\right|\bigg|\sum^{N+1}_{s=2}2^{-s}\sum_{\frac{7}{8}2^{s}<|k|<\frac{9}{8}2^{s}}1\bigg|\gtrsim N\left|\Phi_0(x)\right|.
 \end{equation*}
This immediately leads to \eqref{eq6.1}.

We next prove that if $1<p<\infty$ then
\begin{equation}
\|f_N\|_p\lesssim N^{1/p}.  \label{eq6.8}
\end{equation}
Notice that
 \begin{align*}
 \left|f_N(x)\right|
 &=\left|\Phi_0(x)\right|\left|\sum^{N+1}_{s=2}2^{-s(1-\frac 1p)}\sum_{\frac{7}{8}2^{s}<|k|<\frac{9}{8}2^{s}}e^{i4kx}\right|\\
 &\lesssim \left|\Phi_0(x)\right|\sum^{N+1}_{s=2}2^{-s(1-\frac 1p)}\left(\left|D_{\frac{7}{8}2^{s}}(4x)\right|+\left|D_{\frac{9}{8}2^{s}}(4x)\right|+1\right),
 \end{align*}
where $D_A$ is the Dirichlet kernel,
\begin{equation*}
D_A(t)=\sum\limits_{|k|\leq A}e^{ikt}.
\end{equation*}
It is well-known that
\begin{equation*}
\left|D_A(t)\right|\leq \min\{2A+1, |\sin (t/2)|^{-1}\}
\end{equation*}
(see for example \cite[P. 3]{D01}). Thus
 \begin{equation}
 \left|f_N(x)\right|\lesssim \left|\Phi_0(x)\right|\sum^{N+1}_{s=2}2^{-s(1-\frac 1p)}\min\{2^{s},|\sin (2x)|^{-1}\}.\label{eq6.2}
 \end{equation}
We claim that for any $k\in\mathbb{Z}$
 \begin{equation}\label{eq6.7}
 I_{k,p,N}:=\!\int^{\frac{\pi}{4}(k+1)}_{\frac{\pi}{4}k} \!\!\left(\sum^{N+1}_{s=2}2^{-s(1-\frac 1p)}\min\{2^{s},|\sin (2x)|^{-1}\}\right)^{\!\!p} \textrm{d}x \lesssim_p N.
 \end{equation}
Since $\Phi_0\in\mathscr{S}$, by \eqref{eq6.2} and \eqref{eq6.7} we get
 \begin{equation*}
  \|f_N\|_p\lesssim_p \left( \sum_{k\in\mathbb{Z}} I_{k,p,N}(1+|k|)^{-2p} \right)^{1/p}\lesssim_p N^{1/p},
 \end{equation*}
namely \eqref{eq6.8}.

It follows from \eqref{eq6.1} and \eqref{eq6.8} that $T_a$ is unbounded on $L^p$ if $1<p<\infty$.

It remains to prove \eqref{eq6.7}.  We may assume that $k=0$ since all other cases can be reduced to this case due to the periodicity of the integrand. If $0\leq l\leq N-2$ and  $2^{-l-1}\leq y\leq 2^{-l}$ then
 \begin{align}
\sum^{N+1}_{s=2}2^{-s(1-\frac 1p)}\min\{2^{s},4y^{-1}\} &\leq \sum^{l+2}_{s=2}2^{-s(1-\frac 1p)}2^{s}+\sum^{N+1}_{s=l+3}2^{-s(1-\frac 1p)}4y^{-1}\nonumber\\
 &\lesssim  2^{l/p},\label{eq6.3}
 \end{align}
in which we have used the assumption $p>1$. On the other hand, if $0<y\leq 2^{-N+1}$ then
 \begin{equation}
\sum^{N+1}_{s=2}2^{-s(1-\frac 1p)}\min\{2^{s},4y^{-1}\}\lesssim \sum^{N+1}_{s=2}2^{-s(1-\frac 1p)}2^{s}\lesssim 2^{N/p}.\label{eq6.4}
 \end{equation}
By changing variables and using \eqref{eq6.3} and \eqref{eq6.4},  we get
 \begin{align*}
  I_{0,p,N}&\leq\int^1_0  \left(\sum^{N+1}_{s=2}2^{-s(1-\frac 1p)}\min\{2^{s},y^{-1}\}\right)^{p} \textrm{d}y\\
  &\leq \left(\sum_{l=0}^{N-2}+\sum_{l=N-1}^{\infty}\right) \int_{2^{-l-1}}^{2^{-l}} \left(\sum^{N+1}_{s=2}2^{-s(1-\frac 1p)}\min\{2^{s},4y^{-1}\}\right)^{p} \textrm{d}y\\
 &\lesssim  \sum_{l=0}^{N-2}\int^{2^{-l}}_{2^{-l-1}}2^l \,\textrm{d}y+\sum^{\infty}_{l=N-1}\int^{2^{-l}}_{2^{-l-1}}2^N \,\textrm{d}y\lesssim N,
 \end{align*}
as desired. \qed


%
%
%
%



\begin{thebibliography}{99}


\bibitem{AH90} 
 J. \'{A}lvarez and J. Hounie,
 \emph{Estimates for the kernel and continuity properties of pseudo-differential operators},
 Ark. Mat.,
 28 (1990), no. 1, 1–-22.

\bibitem{Ching} 
C. H. Ching,
\emph{Pseudo-differential operators with nonregular symbols},
J. Differential Equations,
11 (1972), 436-–447.



 \bibitem{D01} 
 J. Duoandikoetxea,
 \emph{Fourier analysis.}
 Translated and revised from the 1995 Spanish original by David Cruz-Uribe. Graduate Studies in Mathematics, 29. American Mathematical Society, Providence, RI, 2001.



 \bibitem{G09} 
 L. Grafakos,
 \emph{Modern Fourier analysis. Second edition.}
 Graduate Texts in Mathematics,  250. Springer, New York, 2009.



 \bibitem{H65} 
 L. H\"{o}rmander,
 \emph{Pseudo-differential operators},
 Comm. Pure Appl. Math.,
 18 (1965), 501–-517.


 \bibitem{H66} 
 L. H\"{o}rmander,
 \emph{Pseudo-differential operators and hypoelliptic equations},
 Singular integrals (Proc. Sympos. Pure Math., Vol. X, Chicago, Ill., 1966), 138–-183. Amer. Math. Soc., Providence, R.I., 1967.




 \bibitem{H71} 
 L. H{\"o}rmander,
 \emph{On the $L^{2}$ continuity of pseudo-differential operators},
  Comm. Pure Appl. Math.,
 24 (1971),  529--535.

 \bibitem{H86}
 J. Hounie,
 \emph{On the $L^2$-continuity of pseudo-differential operators},
 Comm. Partial Differential Equations,
 11 (1986), no.7, 765--778.

\bibitem{KW07}
        C. E. Kenig and W. Staubach,
        \emph{$\Psi$-pseudodifferential operators and estimates for maximal oscillatory integrals},
        Studia Math.,
        183 (2007), no. 3, 249--258.


 \bibitem{KN65}
 J. J. Kohn and L. Nirenberg,
\emph{An algebra of pseudo-differential operators},
 Comm. Pure Appl. Math.,
 18 (1965), 269–-305.


 \bibitem{R76}
 L. Rodino,
 \emph{On the boundedness of pseudo differential operators in the class $L^{m}_{\rho,1}$},
 Proc. Amer. Math. Soc.,
 58 (1976), 211–-215.



 \bibitem{RLS13} 
 S. Rodr\'{\i}guez-L\'{o}pez and W. Staubach,
 \emph{Estimates for rough Fourier integral and pseudodifferential operators and applications to the boundedness of multilinear operators},
 J. Funct. Anal.,
 264 (2013), no. 10, 2356--2385.


\bibitem{SSS91} 
        A. Seeger, C. D. Sogge, and E. M. Stein,
        \emph{Regularity properties of Fourier integral operators},
        Ann. of Math. (2),
        134 (1991), 231--251.


 \bibitem{S10} 
 A. Stefanov,
 \emph{Pseudodifferential operators with rough symbols},
 J. Fourier Anal. Appl.,
 16 (2010), no. 1, 97–-128.

\bibitem{S93} 
        E. M. Stein,
        \emph{Harmonic analysis: real-variable methods, orthogonality, and oscillatory integrals,
        With the assistance of Timothy S. Murphy},
         Princeton Mathematical Series, 43. Monographs in Harmonic Analysis, III. Princeton University Press, Princeton, NJ, 1993.


 \bibitem{T91} 
 M. E. Taylor,
 \emph{Pseudodifferential operators and nonlinear PDE},
 Progress in Mathematics,  100. Birkh\"{a}user Boston, Inc., Boston, MA, 1991.

\end{thebibliography}
\end{document}